\documentclass[11pt,a4paper,fleqn]{article}

\usepackage{amsmath,amssymb,amsthm,enumerate}
\newtheorem{theorem}{Theorem}

\newtheorem{corollary}{Corollary}

{\theoremstyle{definition}
\newtheorem{definition}{Definition}

\newtheorem*{note*}{Note}
}

\begin{document}
\par\noindent {\LARGE\bf
Alternative (Oriented) singular cochains and the modified cup product\\

\par}

{\vspace{5mm}\par\noindent {\it
Taliya Sahihi~$^\star$
Homayoon Eshraghi~$^\dag$ \par\vspace{2mm}\par}

{\vspace{2mm}\par\noindent \it
School of Physics,
Iran University of Science and Technology (IUST), Tehran, Iran}

{\vspace{5mm}\par\noindent {\it
Ali Taghavi~$^\ddag$ \par\vspace{2mm}\par}

{\vspace{2mm}\par\noindent \it
Faculty  of  Mathematics  and Computer Science,
Damghan University, Damghan, Iran}

{\vspace{2mm}\par\noindent {$\phantom{\dag}$~\rm $^\star$\,\, }{\it
taliya\_sahihi@alumni.iust.ac.ir } \par}

{\vspace{2mm}\par\noindent {$\phantom{\dag}$~\rm $^\dag$\,\, }{\it
eshraghi@iust.ac.ir } \par}

{\vspace{2mm}\par\noindent {$\phantom{\dag}$~\rm $^\ddag$\,\, }{\it
taghavi@du.ac.ir } \par}

{\vspace{8mm}\par\noindent\hspace*{10mm}\parbox{140mm}{\small

A special subcomplex of the singular chain complex for a topological space, historically
called oriented singular chain complex is used here with the new name ``alternative" singular
chain complex. It was already known that this subcomplex and so its dual complex are chain homotopy
equivalent to singular chains and cochains respectively and thus have the same homology and
cohomology. Here, in addition to revisiting some aspects of this subcomplex, it is shown that
alternative singular cochains (dual of alternative singular chains) with coefficients in rational
or real numbers are indeed summands of singular cochains through a natural splitting. It is shown
that this natural splitting also hold for cohomologies: At any order, the singular cohomology splits
into the alternative cohomology and another summand which is zero if the considered topological space
is compact. Also in this case similar to the wedge product for differential forms, a modified cup 
product can be defined with the same algebraic properties as in the wedge product in differential forms. 
This provides an idea to investigate some topological and structure-free aspects of nonlinear global 
differential equations on manifolds.

}\par\vspace{6mm}}

\vspace{5mm} \hspace{2mm} Key words: Alternative Chains and Cochains, Alternative Cup Product.

\hspace{2mm} Mathematics Subject Classification (2010): 55N10.

\section{Introduction}

In spite of a huge and almost full developments in homology and cohomology theory in geometry
and algebra, they look not significantly effective in nonlinear features on manifolds and topological spaces.
For example in a differentiable manifold one may ask about solutions of a nonlinear differential equation
of differential forms but another question is: Are these solutions independent from the smooth structure
considered for the manifold? Of course the sense of independence should be explained precisely.
In other words, can these solutions suggest some topological invariants common for all smooth structures?
One may pass to the topological consideration by formally replacing differential forms by singular cochains,
the differentiation operator by the coboundary  operator and the wedge product by the cup product.
The difficulty is that not only the space of singular cochains are totally different from differential forms but also
the cup product does not allow a signed commutation held in the wedge product. In the present note a special
well known subcomplex of singular cochains which we call it here as {\it alternative singular cochains} is
used although historically it was known as "oriented" singular cochains \cite{Barr}. It will be seen that when
the coefficients of singular cochains are rational or real numbers, they naturally split such that alternative
cochains are one of their summands. Also similar to the wedge product in differential forms, one can define
a modified cup product with exactly the same algebraic properties as in the wedge product. This perhaps provides
a step towards determining structure-free features of global nonlinear equations on smooth manifolds. In future reports
by the authors this result will be used to investigate the solutions of a simple nonlinear differential equation
originated in physics.

Historically differential forms were discovered (long before homology and cohomology) in the works of Clairaut
on differential equations in $18$\textsuperscript{th} century \cite{Clairaut1}, \cite{Clairaut2}. In the last decade of $19$\textsuperscript{th} century
 \'{E}lie Cartan represented the differential forms in the language of algebra \cite{Cartan}. After that Poincar\'{e}
proved his lemma about the exactness of closed differential forms on contractible subsets of the Euclidean
space \cite{Katz}. On the other hand at the beginning of $20$\textsuperscript{th} century Poincar\'{e} discovered homological
computation for Betti numbers and torsion coefficients on simplicial and cell structures without knowing
about homology group \cite{Hatcher}.

The study of cohomology from a differentiable viewpoint was almost stopped until 1931 when de Rham found
the relation between homology and differential forms and proved his famous theorem.

However in the middle of 1920s Alexander and Lefschetz found the intersection theory which explains the intersection
of the images of continuous and perhaps smooth functions into the manifold which represent homological elements
\cite{Dieudonne}. Further investigations on differentiable manifolds like the Hodge decomposition theorem appeared later.

As mentioned above, topological viewpoints of homology and cohomology were started by Poincar\'{e} and developed
by others but until the 1920s only simplicial homology and cohomology were known. Then within the next 20 years
many improvements were made by many mathematicians such as Alexander, Lefschetz, Pontryagin, Kolmogorov, Steenrod,
 \v{ Cech} and Whitney and finally Eilenberg constructed the precise singular homology and cohomology theory
\cite{Dieudonne}-\cite{Spanier}.

It is seen that one of the main results relating the differential viewpoint to pure topological one is the
de Rham theorem establishing an isomorphism between the de Rham and singular cohomologies.
Indeed, this theorem relates the solution of differential equation $d\alpha=0$ on a smooth manifold to the solution
of the algebraic equation $\delta\alpha=0$ where in the first equation $\alpha$ is a smooth differential form and
in the second equation  $\alpha$ is a singular cochain and $\delta$ is the cobouandary operator. As mentioned
before passing to nonlinear equations meets many difficulties for example we may ask what is the relation between
the solutions of equation $\alpha \wedge d\alpha=0$  on a smooth manifold and the solutions of $\alpha \smile\delta\alpha=0$
in terms of cochains on the same manifolds?

Because the cup product on arbitrary cochains does not allow the signed commutativity governing the wedge product,
it looks no significant relation between these two equations. However if we restrict the cochain group at any order
to its subgroup, alternative cochains, one can easily modify the cup product to be very similar to the wedge
product. In the near future the authors will give a detailed comparison between the equations $\alpha \wedge d\alpha=0$
and $\alpha \smile\delta\alpha=0$ by the use of alternative cochains.

As mentioned in the abstract, alternative chains and cochains are indeed subcomplexes of singular chains and cochains
historically known as oriented singular chains and cochains respectively and thus admit alternative homology and cohomology.
Although at first the oriented (alternative) singular chain complex was used by Lefschetz, fainally it was shown by Barr in 1995
that the chain complex of singular chains is chain homotopy equivalent with its subcomplex oriented (alternative)
singular chains and thus alternative homology and cohomology are respectively isomorphic with singular ones \cite{Barr}.
These isomorphisms are very good news for our purposes in cochains since we can only work with alternative cochains without
loosing any thing about their cohomologies.

In the next section (Sec. 2) alternative singular cochains are defined for arbitrary groups or fields of coefficients. Then especially
when the coefficients are in rational numbers $\mathbb{Q}$  or real numbers $\mathbb{R}$, a homomorphism called {\it alternative maker
operator} denoted by $A$ is defined which splits the group of singular cochains to two subgroups at any order. Using this
operator enables us to modify the cup product on alternative cochains to make it similar to the wedge product.

Section 3 performs two main tasks. The first is the existence of the alternative cohomology at any positive order
and with any coefficient group, that is, alternative cochains indeed form a subcomplex of singular cochains. This
has been known since 1940s which was also reviewed in Ref. \cite{Barr}. However the approach presented here is rather
different yielding a more direct proof. The second is that the alternative maker operator $A$ commutes with coboundary
operator $\delta$ and so induces the same homomorphism on cohomology groups and the same splitting of the cohomology groups.
However due to the isomorphism between the singular and alternative cohomologies, it will be seen that if the considered
topological space is compact, $A$ is indeed an isomorphism between the alternative and singular cohomologies. The situation
for noncompact spaces is not yet known for the authors. 

In Sec. 4 alternative chains are defined leading to the existence of the alternative homology as known from 1940s but
again our approach is more direct and easier to follow. Although the chain homotopy equivalence established by Barr in
Ref. \cite{Barr} of course proves the isomorphism between alternative and singular homologies, but here a quite different
proof using simplicial homology is employed which is very similar to the well known proof for singular and simplicial
homologies. Indeed, it is shown that for a simplicial complex this homology is isomorphic with simplicial one and
consequently isomorphic with singular homology. Since any topological space is weakly  homotopy equivalent to a simplicial
complex, these homologies are isomorphic for any space. However, our approach has a disadvantage that it does not imply
an isomorphism between cohomologies because alternative chain groups are not free and have torsion summands. The last
section (Sec. 5) summarizes the results together with some remarks and suggestions for further investigations.

\hspace{0.25cm}
\section{Alternative (Oriented) singular cochains}

Simplices as smallest convex subsets of Euclidean spaces are the best tools to construct homology and cohomology groups.
Consider $n+1$ points $v_{0}, v_{1},\ldots,v_{n}$ in $\mathbb{R}^{n}$ such that all vectors $v_{1}-v_{0},v_{2}-v_{0},\ldots,v_{n}-v_{0}$
are linearly independent. An $n$-simplex $\Delta^{n}$ is defined as follows

\begin{equation}\label{simplex}
 \Delta^{n}=\{\sum_{j=0}^{n} t_{j}v_{j}\mid \sum_{j=0}^{n} t_{j}=1,\forall j: t_{j}\geq0\} \subset\mathbb{R}^{n}.
\end{equation}

The points $v_{0}, v_{1},\ldots,v_{n}$ are the vertices of the simplex. When the order of the vertices is important, the simplex
is called an ordered $n$-simplex and denoted by $[v_{0},\ldots,v_{n}]$.
Each $n$-simplex has $n+1$ faces $F_{i}(\Delta^{n})$ ($i=0,\dots,n$) as

\begin{equation}\label{face}
F_{i}(\Delta^{n})=\{\sum_{j=0}^{n} t_{j}v_{j}\mid\sum_{j=0}^{n} t_{j}=1,t_{i}\equiv 0,\forall j\neq i: t_{j}\geq0 \},\\
(i=0,1,\ldots,n).
\end{equation}

For an ordered simplex the $i$\textsuperscript{th} face $F_{i}$ is denoted by  $[v_{0},\ldots,\hat{v}_{i},\ldots,v_{n}]$
meaning that the vertex $v_{i}$ is omitted. The order of each face is obviously induced by the simplex order.
Because later the faces will be used, let us explicitly denote $F_{i}$ by

\begin{equation}\label{reorder}
F_{i}=[v_{0},\ldots,\hat{v}_{i},\ldots,v_{n}]=[w_{0},\ldots,w_{n-1}], \\ w_{j}=\begin{cases}
v_{j}: & 0\leq j<i\\
v_{j+1}: & i\leq j\leq n-1
\end{cases}.
\end{equation}

In this note everywhere only ordered simplices are considered and so any permutation of the vertices induces
a homeomorphism on the simplex. To formulate this, assume $S_{k}$ be the group of permutations of $k$
things i.e., the set of all bijections on  $\{1,\ldots,k\}$ and the group operation is the combination of bijections.

\begin{definition}
For an ordered $n$-simplex $\Delta^{n}=[v_{0},\ldots,v_{n}]$, each permutation  $s \in S_{n+1}$ induces a permutation transformation
$\overline{s}: \Delta^{n} \rightarrow \Delta^{n}$ defined as

\begin{equation}\nonumber
x=\sum_{j=0}^{n} t_{j}v_{j}\longmapsto \bar{s}(x)=\sum_{j=0}^{n} t_{j}v_{s(j)}=\sum_{j=0}^{n} t_{s^{-1}(j)}v_{j}.
\end{equation}
\end{definition}
Clearly, $(\overline{s})^{-1}=\overline{(s^{-1})}$ and so $\overline{s}$ is a homeomorphism.

Let us denote the free abelian group generated by all continuous maps $\sigma:\Delta^{n}\rightarrow X$
by $C_{n}(X)$ for the space $X$ which is also called the group of $n$-chains such as $\sum_{i}^{} n_{i}\sigma_{i}$.
The permutation transformation $\overline{s}$ induces a homomorphism    $\overline{s}^{\ast}:C_{n}(X)\rightarrow C_{n}(X)$  through the
 pull back of each generator $\sigma:\Delta^{n}\rightarrow X$ as

\begin{equation}\label{perm2}
\overline{s}^{\ast}\sigma(\sum_{j}^{} t_{j}v_{j})=\sigma\circ\bar{s}(\sum_{j}^{} t_{j}v_{j})=\sigma(\sum_{j}^{} t_{j}v_{s(j)}).
\end{equation}

As usual the boundary operator (homomorphism) $\partial:C_{n}(X)\rightarrow C_{n-1}(X)$  for each generator $\sigma$
is defined as $\partial\sigma=\sum_{j}^{}(-1)^{j}\sigma \mid_{[v_{0},\ldots,\hat{v}_{j},\ldots,v_{n}]}$. Now we need  to calculate
$\partial\overline{s}^{\ast}\sigma$ explicitly. To do this, let us make a definition.

\begin{definition}
For each permutation $s\in S_{n+1}$ and each $i\in \{0,1,\ldots,n\}$ the permutation  $s_{i}\in S_{n}$ as a bijection on
$\{0,\ldots,n-1\}$ is induced by eliminating $i$ from the set $ \{0,\ldots,n\}$ and $s(i)$ from $ \{s(0),\ldots,s(n)\}$ and restricting
$s$ to these reduced sets, or explicitly for each $j\in \{0,\ldots,n-1\}$ we have

\begin{equation}\nonumber
\begin{split}
s_{i}(j)=\begin{cases}
s(j): & 0\leq j<i,s(j)<s(i),\\
s(j)-1: & 0\leq j<i, s(j)>s(i),\\
s(j+1): & i\leq j\leq n-1,s(j+1)<s(i),\\
s(j+1)-1: & i\leq j\leq n-1,s(j+1)>s(i).
\end{cases}.
\end{split}
\end{equation}

\end{definition}
This definition simply means that: Eliminate $i$ from $ \{0,\ldots,n\}$ and renumerate the remaining set $ \{0,\ldots,\hat{i},\ldots,n\}$ as
$ \{0,\ldots,n-1\}$. Also eliminate $s(i)$ and renumerate the remaining set $ \{0,\ldots,\hat{s}(i),\ldots,n\}$  as $ \{0,\ldots,n-1\}$, the act of $s_{i}$
on these reduced sets is the same as $s$.

The restriction of $\overline{s}^{\ast}\sigma$ to each face $F_{i}$ according to (\ref{perm2}) is defined as:

\begin{equation}\nonumber
\overline{s}^{\ast}\sigma(\sum_{j}^{} t_{j}v_{j})|_{t_{i}=0}=\sigma(\sum_{j}^{} t_{j}v_{s(j)})|_{t_{i}=0}=
\sigma(\sum_{j\neq  i} t_{j}v_{s(j)}),
\end{equation}
 where $\sum_{j\neq i} t_{j}=1$ of course holds. The above equation means that we must change the order of vertices by $s$
 and then eliminate vertex $v_{s(i)}$ and so restrict to the face $F_{s(i)}$ and finally take this face by $\sigma$ to $X$.  Using Eq.~(\ref{reorder}),
 denote the face $F_{s(i)}$ as $[v_{0},\ldots,\hat{v}_{s(i)},\ldots,v_{n}]=[w_{0},\ldots,w_{n-1}]$, then the definition of $s_{i}$ (Definition 2) implies that $\sum_{j\neq i}^{} t_{j}v_{s(j)}=\sum_{j=0}^{n-1} t_{j}w_{s_{i}(j)}$
 and hence

\begin{equation}\nonumber
(\overline{s}^{\ast}\sigma) |_{F_{i}}(\sum_{j=0}^{n} t_{j}v_{j})
=\sigma |_{F_{s(i)}}(\sum_{j=0}^{n-1} t_{j} w_{s_{i}(j)})
=\overline{s}_{i}^{\ast}(\sigma |_{F_{s(i)}})(\sum_{j=0}^{n-1} t_{j} w_{j}).
\end{equation}
Finally we obtain

\begin{equation}\label{5}
\begin{split}
\partial(\overline{s}^{\ast}\sigma)=\sum_{i=0}^{n}(-1)^{i}(\overline{s}^{\ast}\sigma)|_{[v_{0},\ldots,\hat{v}_{i},\ldots,v_{n}]}=
\sum_{i=0}^{n}(-1)^{i} \overline{s}_{i}^{\ast}(\sigma|_{[v_{0},\ldots,\hat{v}_{s(i)},\ldots,v_{n}]}).
\end{split}
\end{equation}
The above equation will be used in the next two sections.

By definition, the group of cochains with coefficients in an arbitrary abelian group $G$ is the dual space
of chains, that is $C_{n}^{\ast}(X)=\mathrm{Hom}(C_{n}(X),G)$ which is usually denoted by $C^{n}(X;G)$. Here for
the briefness we drop $G$ and simply write $C^{n}(X)$. In Eq.~(\ref{perm2}) $\overline{s}^{\ast}$ as a homomorphism
on $C_{n}(X)$ was defined and so $\overline{s}^{\ast\ast}$ is a homomorphism on $C^{n}(X)$ acing as the pull back of
 $\overline{s}^{\ast}$. Each permutation $s\in S_{n+1}$ has a well known sign $\epsilon(s)\in \{-1,+1\}$ according to that the number
 of interchanges is odd or even. Now we are ready to give the definition of alternative cochains.

 \begin{definition}
 A singular cochain $\alpha\in C^{n}(X)$ is  alternative (oriented) if for each $s\in S_{n+1}$ it satisfies

\begin{equation}\nonumber
\overline{s}^{\ast\ast}\alpha =\epsilon(s)\alpha.
\end{equation}
This means that for each generator $\sigma$ of $C_{n}(X)$ we have

\begin{equation}\label{altCOCH}
\alpha(\overline{s}^{\ast}\sigma)=\alpha(\sigma\circ\overline{s})=\epsilon(s)\alpha(\sigma).
\end{equation}
\end{definition}

It is easily observed that alternative cochains are a subgroup of $C^{n}(X)$ which we denote it by $C_{A}^{n}(X)$
and $C_{A}^{0}(X)=C_{}^{0}(X)$.

Alternative cochains would be more efficient if the coefficient group $G$ is the field of rational numbers $\mathbb{Q}$
or real numbers $\mathbb{R}$.
So we assume it in the remaining of this section and define the alternative maker operator $A$ as follows.

\begin{definition}
Alternative maker operator $A\in \mathrm{Hom}(C^{n}_{}(X),C^{n}_{A}(X))$ when the coefficient group is
$\mathbb{Q}$ or  $\mathbb{R}$ is such that for each $\alpha\in C^{n}_{}(X)$:

\begin{equation}\nonumber
A(\alpha):=\frac{1}{(n+1)!}\sum_{s\in S_{n+1}}\epsilon(s)\overline{s}^{\ast\ast}\alpha.
\end{equation}
\end{definition}

It is not difficult to see that $A(\alpha)$ is in $C_{A}^{n}(X)$, in fact for each $\sigma\in C_{n}(X)$ and each $\alpha\in C^{n}(X)$
and each $s' \in S_{n+1}$ we find

\begin{equation}\nonumber
A(\alpha)(\overline{s'}^{\ast}\sigma)=\frac{1}{(n+1)!}\sum_{s\in S_{n+1}}\epsilon(s)\alpha ((\overline{s}^{\ast}\circ \bar{s'}^{\ast})\sigma)=
\frac{1}{(n+1)!}\sum_{s\in S_{n+1}}\epsilon(s)\alpha ((\overline{s'\circ s})^{\ast}\sigma),
\end{equation}
and if we put $s''=s'\circ s$ then $\epsilon(s)= \epsilon(s'^{-1}\circ s'')= \epsilon(s')\epsilon(s'')$ and the right hand side becomes

\begin{equation}\nonumber
\frac{\epsilon(s')}{(n+1)!}\sum_{s''\in S_{n+1}}\epsilon(s'')\alpha ((\overline{s''}^{\ast}\sigma)=\epsilon(s')A(\alpha)(\sigma),
\end{equation}
which proves the assertion. Also it is clear that $A$ becomes identity when restricted to $C_{A}^{n}(X)$. If $\alpha \in C_{A}^{n}(X) $,
from definitions 3 and 4 we have

\begin{equation}\nonumber
A(\alpha)=\frac{1}{(n+1)!}\sum_{s\in S_{n+1}}\epsilon(s)\overline{s}^{\ast\ast}\alpha=\frac{1}{(n+1)!}\sum_{s\in S_{n+1}}\epsilon^{2}(s)\alpha=\alpha.
\end{equation}
Therefore, if $j:C_{A}^{n}(X)\hookrightarrow C_{}^{n}(X)$ be the inclusion, then

\begin{equation}\label{7}
A\circ j=id: C_{A}^{n}(X)\rightarrow{}  C_{A}^{n}(X).
\end{equation}
Then the following result is immediate.

\begin{theorem} For each $n$,

\begin{equation}\nonumber
C^{n}(X)=C^{n}_{A}(X)\oplus Ker A.
\end{equation}
\end{theorem}
Of course for n=0, A is the identity and Ker A=0.

{\it Proof}. It suffices to note to the following short exact sequence

\begin{equation}\label{seq8}
0\xrightarrow{\text{}}Ker A\xrightarrow{\text{j}} C^{n}(X)\xrightarrow{\text{A}}  C^{n}_{A}(X)\xrightarrow{}0.
\end{equation}
The result follows from (\ref{7}) and the splitting lemma \cite{Hatcher}. $\blacksquare$\\

The last definition of this section concerns with the modification of the cup product which is very similar to the wedge
product in differential forms. The common cup product $\alpha\smile\beta$ for $\alpha\in C^p(X)$ and $\beta\in C^q(X)$
is defined by its action on a generator of $C_{p+q}(X)$, say $\sigma:\Delta^{p+q}\xrightarrow{} X$ as
$(\alpha\smile\beta)(\sigma)=\alpha(\sigma |_{[v_{0},\ldots, v_{p}]})\beta(\sigma |_{[v_{p},\ldots, v_{p+q}]})$.
A modified cup product similar to the wedge product can be defined as follows.

\begin{definition} For each $\alpha\in C^p(X)$ and $\beta\in C^q(X)$, the alternative cup product
$\alpha\overset{A}{\smile}\beta$ is defined as

\begin{equation}\nonumber
(\alpha\overset{A}{\smile}\beta):=A(\alpha\smile\beta),
\end{equation}
where $A$ is the alternative maker operator (Definition 4) and $\smile$ is the common cup product.

\end{definition}

It is not difficult to see that this modified cup product admits the associativity and the signed commutation just like
the wedge product. Although the proof is very similar to exterior forms but because of some minor changes, it is
worth to give the proof.

\begin{theorem} For each $\alpha\in C^p_{A}(X)$, $\beta\in C^q_{A}(X)$ and $\gamma\in C^r_{A}(X)$ we have

\begin{equation}\nonumber
\begin{split}
(a)  &&&\beta \overset{A}{\smile} \alpha= (-1)^{pq}\alpha\overset{A}{\smile} \beta,\\
(b)  &&&(\alpha\overset{A}{\smile}\beta)\overset{A}{\smile}\gamma= \alpha\overset{A}{\smile}(\beta\overset{A}{\smile}\gamma).
\end{split}
\end{equation}

 \end{theorem}

 {\it Proof}. To prove (a) let $\sigma:\Delta^{p+q}\xrightarrow{} X$ be a generator of $C_{p+q}(X)$ and put
 $\Delta^{p+q}=[v_{0},\ldots, v_{p+q}]$, then

\begin{equation}\nonumber
\begin{split}
\beta\overset{A}{\smile}\alpha=\frac{1}{(p+q+1)!}\sum_{s\in S_{p+q+1}}\epsilon(s)\beta((\sigma\circ \bar{s}) |_{[v_{0},\ldots, v_{q}]})
\alpha((\sigma\circ \bar{s}) |_{[v_{q},\ldots, v_{p+q}]})\\
=\frac{1}{(p+q+1)!}\sum_{s\in S_{p+q+1}}\epsilon(s)\beta(\sigma |_{[v_{s(0)},\ldots, v_{s(q)}]})
\alpha(\sigma |_{[v_{s(q)},\ldots, v_{s(p+q)}]}).
\end{split}
\end{equation}

It is easy to find that a permutation which transforms $(s(0), \dots, s(q))$ to  $(s(q), s(0), \dots, s(q-1))$ has the sign $(-1)^{q}$
and similarly $(-1)^{p}$ is the sign of a permutation changing $(s(q), \dots, s(p+q))$ to $(s(q+1),\dots, s(p+q) ,s(q))$. Since
$\alpha$ and $\beta$ are alternative, the right hand side of the above equation becomes

\begin{equation}\nonumber
\frac{(-1)^{p+q}}{(p+q+1)!}\sum_{s\in S_{p+q+1}}\epsilon(s)\alpha(\sigma |_{[v_{s(q+1)}, \dots, v_{s(p+q)}, v_{s(q)}]})
\beta(\sigma |_{[v_{s(q)}, v_{s(0)},\dots, v_{s(q-1)}]}).
\end{equation}

Now let $s'\in S_{p+q+1}$ be a permutation such that $(s'(0),\dots, s'(p+q))=(s(q+1),\dots, s(q+p), s(q), s(0),\dots, s(q-1))$, then
it is observed that $\epsilon(s)=(-1)^{pq+p+q}\epsilon(s')$ and so

\begin{equation}\nonumber
\beta \overset{A}{\smile} \alpha=\frac{(-1)^{pq}}{(p+q+1)!}\sum_{s'\in S_{p+q+1}}\epsilon(s')\alpha(\sigma |_{[v_{s'(0)},\ldots, v_{s'(p)}]})
\beta(\sigma |_{[v_{s'(p)}, \ldots, v_{s'(p+q)}]})=(-1)^{pq}\alpha\overset{A}{\smile}\beta .
\end{equation}

Now the proof of (b) follows if we transform both sides of (b) to become equal. Let begin with the left hand side

\begin{equation}\nonumber
\begin{split}
(\alpha\overset{A}{\smile}\beta)\overset{A}{\smile}\gamma=\frac{1}{(p+q+r+1)!}\sum_{s\in S_{p+q+r+1}}\epsilon(s)(\alpha\overset{A}{\smile}\beta)
(\sigma |_{[v_{s(0)},\dots,  v_{s(p+q)}]})\gamma(\sigma |_{[v_{s(p+q)}, \dots, v_{s(p+q+r)}]}).
\end{split}
\end{equation}
Assume $u \in S_{p+q+1}$ to be represented as an element of $S_{p+q+r+1}$ such that $u(i)=i$ for $i>p+q$ and therefore $u$
interchanges all numbers $0\leq i \leq p+q$ with each other. Thus indeed $u$ belongs to a subgroup of $S_{p+q+r+1}$ which
we again denote it by $S_{p+q+1}$. Since $\alpha\overset{A}{\smile}\beta$ is an alternative cochain, any permutation induced
by $u$ inserts the sign $\epsilon(u)$ which together with the fact that $\epsilon(u\circ s)=\epsilon(u)\epsilon(s)$ leads to

\begin{equation}\nonumber
\begin{split}
(\alpha\overset{A}{\smile}\beta)\overset{A}{\smile}\gamma=\frac{1}{(p+q+r+1)!(p+q+1)!}\sum_{s\in S_{p+q+r+1}}\sum_{u\in S_{p+q+1}}
\epsilon(u\circ s)\\
\alpha(\sigma |_{[v_{u\circ s(0)},\ldots,  v_{u\circ s(p)}]})
\beta(\sigma |_{[v_{u\circ s(p)},\ldots,  v_{u\circ s(p+q)}]})
\gamma(\sigma |_{[v_{u\circ s(p+q)}, \ldots, v_{u\circ s(p+q+r)}]}).
\end{split}
\end{equation}
Now fix an element $s'\in S_{p+q+r+1}$ and ask
how many pairs $(s,u)$ exist such that $u\circ s=s'$? Of course for any $u\in S_{p+q+1}$ there is one and only one
$s\in S_{p+q+r+1}$ such that $u\circ s=s'$. On the other hand for any fixed pair $(s,u)$  there is one and only one $s'$.
The number of elements $u$ is $(p+q+1)!$ and thus for any fixed $s'$, there are exactly $(p+q+1)!$ pairs $(s,u)$ such that
$u\circ s=s'$. Hence the double summation is reduced to a single summation:

\begin{equation}\nonumber
\sum_{s\in S_{p+q+r+1}}\sum_{u\in S_{p+q+1}}=(p+q+1)!\sum_{s'\in S_{p+q+r+1}}.
\end{equation}
Therefore

\begin{equation}\nonumber
\begin{split}
(\alpha\overset{A}{\smile}\beta)\overset{A}{\smile}\gamma=\frac{1}{(p+q+r+1)!}\sum_{s'\in S_{p+q+r+1}}
\epsilon(s')\alpha(\sigma |_{[v_{s'(0)},\ldots,  v_{s'(p)}]})\beta(\sigma |_{[v_{s'(p)},\ldots,  v_{s'(p+q)}]})\\
\gamma(\sigma |_{[v_{s'(p+q)}, \ldots, v_{s'(p+q+r)}]}).
\end{split}
\end{equation}
A similar procedure transforms the right hand side of (b) to the above relation and the proof is achieved. $\blacksquare$\\

In order to make a significant comparison between nonlinear differential equations and nonlinear cochain equations,
it is necessary that differential forms and cochains obey the same algebraic properties. Theorem 2 establishes two
of these properties and the last one dealing with the coboundary operator $\delta$ is proved in the next section.

The last important point in this section is that the splitting of theorem 1 is canonical or natural in the sense that for any
continuous map $f:X\xrightarrow{}Y$, $f^{\ast}$ sends each element of $C^{n}_{A}(Y)$ to $C^{n}_{A}(X)$ and also
$f^{\ast}A=Af^{\ast}$. This naturality clearly extends to the solutions of any nonlinear cochain equation provided that $A$
also commutes with $\delta$ as seen in the next section.

\hspace{0.25cm}
\section{Alternative (Oriented) singular cohomology}

As mentioned before, the alternative singular cohomology was known from 1940s although its isomorphism with the singular
cohomology was proved in 1995 \cite{Barr}. However, in the following theorem a more direct and straightforward demonstration is
presented for the existence of the alternative (oriented) singular cohomology.

\begin{theorem} For a topological space $X$ and with coefficients in any abelian group, the following sequence is a
chain complex in which $\delta$ is the coboundary operator.

\begin{equation}\nonumber
0\xrightarrow{\text{}} C^{0}_{A}(X)\xrightarrow{\text{$\delta_{}$}} C^{1}_{A}(X)\xrightarrow{\text{$\delta_{}$}}\ldots \xrightarrow{\text{$\delta_{}$}} C^{n-1}_{A}(X)\xrightarrow{\text{$\delta_{}$}} C^{n}_{A}(X)\xrightarrow{\text{$\delta_{}$}} C^{n+1}_{A}(X)\xrightarrow{\text{$\delta_{}$}} \ldots.
\end{equation}
Therefore the alternative cohomology $\mathrm{H}_{A}^{n}(X)=Ker \delta / Im \delta$  is defined.
\end{theorem}

{\it Proof}. It is sufficient to prove that $\delta$ takes  $C^{n}_{A}(X)$ to $C^{n+1}_{A}(X)$ for $n\geq0$. First we note
that $C^{0}_{A}(X)=C^{0}_{}(X)$ and if $\alpha \in C^{0}_{}(X)$ and $\sigma : \Delta ^{1}=[v_{0},v_{1}]\xrightarrow{} X$
be a generator of $C_1(X)$ then $\delta\alpha(\sigma)=\alpha(\sigma|_{[v_0]})-\alpha(\sigma|_{[v_1]})$. Assume $s\in S_{2}$ is the only nontrivial
element of this group such that $s(0)=1$ and $s(1)=0$. Then

\begin{equation}\nonumber
\begin{split}
(\overline{s}^{\ast\ast}\delta\alpha)(\sigma)=\delta\alpha(\overline{s}^{\ast}\sigma)=\delta\alpha(\sigma\circ s)=\delta\alpha(\sigma|_{[v_1,v_0]})
=\alpha(\sigma|_{[v_1]})-\alpha(\sigma|_{[v_0]})\\
=-\delta\alpha(\sigma)=\epsilon(s)\delta\alpha(\sigma).
\end{split}
\end{equation}

Now assume $\alpha\in C^{n-1}_{A}(X)$ and $\sigma:\Delta^{n}\xrightarrow{} X$ for $n\geq 2$. For any $s\in S_{n+1}$

\begin{equation}\nonumber
\begin{split}
(\overline{s}^{\ast\ast}\delta\alpha)(\sigma)=\delta\alpha(\overline{s}^{\ast}\sigma)=\alpha(\partial\overline{s}^{\ast}\sigma).
\end{split}
\end{equation}
Applying Eq.~(\ref{5}) yields

\begin{equation}\nonumber
\begin{split}
(\overline{s}^{\ast\ast}\delta\alpha)(\sigma)=\sum_{i=0}^{n}(-1)^{i}\alpha(\overline{s}_{i}^{\ast}(\sigma|_{[v_{0},\ldots,\hat{v}_{s(i)},\ldots,v_{n}]})),
\end{split}
\end{equation}
where $s_{i}$ is given by Definition 2. Since $\alpha$ is alternative, the right hand side becomes

\begin{equation}\nonumber
\begin{split}
\sum_{i=0}^{n}(-1)^{i}\epsilon(s_i)\alpha(\sigma|_{[v_{0},\ldots,\hat{v}_{s(i)},\ldots,v_{n}]}).
\end{split}
\end{equation}

Now we must drive the relation between $\epsilon (s_i)$ and $\epsilon (s)$. Indeed, $s$ puts $s(i)$ in the $i^{\text{\tiny th}}$
place and the remaining permutations are the same as for $s_i$ . Taking $s(i)$ from its original place to the $i^{\text{\tiny th}}$
place  contributes the sign $(-1)^{s(i)-i}$ in $\epsilon(s)$ and so
\begin{equation}\label{sigma}
\epsilon(s)=(-1)^{i-s(i)}\epsilon(s_{i}).
\end{equation}
Therefore we obtain
\begin{equation}\nonumber
(\overline{s}^{\ast\ast}\delta\alpha)(\sigma)=\epsilon(s)\sum_{i=0}^{n}(-1)^{s(i)}\alpha(\sigma|_{[v_{0},\ldots ,\hat{v}_{s(i)},\ldots,v_{n}]})=\epsilon(s)\delta\alpha(\sigma).
\end{equation}
Thus $\delta\alpha\in C_{A}^{n}(X)$. $\blacksquare$
It is also seen that for $n=0$, the alternative cohomology is the same as usual singular cohomology
\begin{equation}\label{iso}
\mathrm{H}_{A}^{0}(X)=\mathrm{H}_{}^{0}(X).
\end{equation}

If $f:X\xrightarrow{}Y$ be a continuous map between spaces $X$ and $Y$,  since $f^{\ast}$ commutes with $\delta$ and takes
$C_{A}^{n}(Y)$ into $C_{A}^{n}(X)$, $f^{\ast}:\mathrm{H}_{A}^{n}(Y)\xrightarrow{} \mathrm{H}_{A}^{n}(X)$ is defined. Also this cohomology
is defined for relative cohomology such as $\mathrm{H}_{A}^{n}(X,B)$ for $B \subset X$. Long exact sequences for pairs $(X,B)$ or for triple
$(X,B,C)$ with $C\subset B\subset X$ are all valid for the alternative cohomology as well. However the equality $f^{\ast}=g^{\ast}$
when $f$ and $g$ are homotopic is established in the next section after introducing the alternative homology.

To see the commutation of $A$ with $\delta$ and its subsequent results let the group of coefficients of cochains be $\mathbb{Q}$
or $\mathbb{R}$.

\begin{theorem} When the coefficients are in $\mathbb{R}$ or $\mathbb{Q}$, the alternative maker operator and the coboundary
operator commute on $C^n(X)$ for all $n$:
\begin{equation}\nonumber
A\delta=\delta A.
\end{equation}
\end{theorem}
{\it Proof}. For each $\alpha\in C^{n-1}(X)$ and $\sigma:\Delta^n\xrightarrow{} X$ a generator of $C_n(X)$ we have

\begin{equation}\nonumber
\begin{split}
\delta(A(\alpha))(\sigma)=A(\alpha)(\partial\sigma)=\frac{1}{n!}\sum_{s'\in S_n}\epsilon(s')\alpha(\overline{s'}^{\ast}(\partial\sigma))\\
=\frac{1}{(n+1)!}\sum_{s'\in S_n}(n+1)\epsilon(s')\sum_{j=0}^{n}(-1)^{j}\alpha(\bar{s'}^{\ast}(\sigma |_{[v_{0},\dots, \hat{v}_{j},\dots, v_{n}]})).
\end{split}
\end{equation}

For each fixed $s'\in S_n$ and fixed $j\in \{0,\dots, n\}$ in this summation we can define $(n+1)$ permutations
$s^{(0)}(j,s'),\ldots, s^{(n)}(j,s')$ in $S_{n+1}$ such that $s^{(i)}_{i}(j,s')=s'$ for each $i\in \{0,\dots, n\}$. The
lower index $i$ is given by Definition 2. This means that for each $i$, $(s^{(i)}_{}(j,s'))(i)=j$, that is, it takes $i$ to
the fixed $j$ and so from Eq.~(\ref{sigma}) it follows that

\begin{equation}\nonumber
\epsilon(s')=(-1)^{i-j}\epsilon(s^{(i)}(j,s')),
\end{equation}
and therefore

\begin{equation}\nonumber
(n+1)\epsilon(s')=\underbrace{\epsilon(s')+\dots +\epsilon(s')}_{\text{ (n+1)-times }}=\sum_{i=0}^{n}(-1)^{i-j}\epsilon(s^{(i)}(j,s')).
\end{equation}
Consequently

\begin{equation}\nonumber
\delta(A(\alpha))(\sigma)=\frac{1}{(n+1)!}\sum_{s'\in S_{n}}\sum_{j=0}^{n}\sum_{i=0}^{n}(-1)^{i}\epsilon(s^{(i)}(j,s'))\alpha(\overline{s}^{(i)\ast}_{i}(j,s')
(\sigma |_{[v_{0},\dots,\hat{v}_{j}\dots, v_{n}]})),
\end{equation}
and keep in mind that $(s^{(i)}(j,s'))(i)=j$. For any fixed permutation $s\in S_{n+1}$ how many $i$ , $j$ and $s'$ do exist such that
$s^{(i)}(j,s')=s$? Indeed for each $i$, $j$  is determined by $s(i)=j$ and $s'$ is also determined by $s_{i}=s'$. In  other words, for
each $i$, there is a one to one correspondence between the pairs $(j,s')$ and $s\in S_{n+1}$. Therefore the double sum
$\sum_{s'\in S_{n}}\sum_{j=0}^{n} $ is replaced with the single sum $ \sum_{s \in S_{n+1}}$ and hence

\begin{equation}\nonumber
\delta(A(\alpha))(\sigma)=\frac{1}{(n+1)!}\sum_{s\in S_{n+1}}\epsilon(s)\sum_{i=0}^{n}(-1)^{i}\alpha(\overline{s}^{\ast}_{i}(\sigma |_{[v_{0},\dots,\hat{v}_{s(i)}\dots, v_{n}]})).
\end{equation}

 Finally using Eq.~(\ref{5}), we obtain

\begin{equation}\nonumber
\delta(A(\alpha))(\sigma)=\frac{1}{(n+1)!}\sum_{s\in S_{n+1}}\epsilon(s)\alpha(\partial(\overline{s}^{\ast}\sigma))
=\frac{1}{(n+1)!}\sum_{s\in S_{n+1}}\epsilon(s)\overline{s}^{\ast\ast}\delta\alpha(\sigma)=A(\delta\alpha)(\sigma).
\end{equation}  $\blacksquare$\\

This theorem gives two important results. The first is the last property of alternative cup product.

\begin{corollary} When the coefficients are in $\mathbb{R}$ or $\mathbb{Q}$, alternative cup product satisfies

\begin{equation}\nonumber
\delta(\alpha\overset{A}{\smile}\beta)=\delta\alpha\overset{A}{\smile}\beta+(-1)^{p}\alpha\overset{A}{\smile}\delta\beta,
\end{equation}
where $\alpha\in C^{p}_{A}(X)$ and $\beta$ is an arbitrary alternative cochain.
 \end{corollary}

{\it Proof}. Since the ordinary cup product satisfies this property,

\begin{equation}\nonumber
\delta(\alpha\overset{A}{\smile}\beta)=\delta A(\alpha\smile\beta)=A(\delta (\alpha\smile\beta))=A(\delta\alpha\smile\beta+(-1)^{p}\alpha\overset{}{\smile}\delta\beta)=
\delta\alpha\overset{A}{\smile}\beta+(-1)^{p}\alpha\overset{A}{\smile}\delta\beta.
\end{equation} $\blacksquare$\\
Therefore a comparison between the differential forms and the alternative cochains seems to have sense.
The second result is the splitting of the cohomology groups.

\begin{corollary} When the coefficients are in $\mathbb{R}$ or $\mathbb{Q}$, the alternative maker operator
$A$ implies the same splitting as in Theorem 1 for cohomology groups:
\begin{equation}\nonumber
\mathrm{H}^{n}(X)=\mathrm{H}_{A}^{n}(X)\oplus Ker A.
\end{equation}
Moreover, if the singular cohomology is finitely generated, for example if the space $X$ is compact, then
$A$ becomes an isomorphism.
\end{corollary}
{\it Proof}. The last theorem $(\delta A=A\delta)$ implies that $A$ is defined also on the cohomology groups,
$A:\mathrm{H}^{n}(X)\xrightarrow{} \mathrm{H}_{A}^{n}(X)$ and so Eq.~(\ref{7}) and the short exact sequence (\ref{seq8})
also hold and thus the splitting of Theorem 1 also holds for cohomology level. On the other hand due to Ref. \cite{Barr}
$\mathrm{H}^{n}(X)$ is isomorphic to $\mathrm{H}_{A}^{n}(X)$ and if $\mathrm{H}^{n}(X)$ is finitely generated,
the fundamental theorem for finitely generated abelian groups leads to $Ker A=0$ and thus $A$ (already seen to be
surjective) becomes an isomorphism between the singular and alternative cohomologies. However it should be
mentioned that here the cohomology groups are indeed $\mathbb{R}$-modules or $\mathbb{Q}$-modules but the result
is still kept. $\blacksquare$\\
It is noticeable that if the singular cohomology is not finitely generated for some non compact spaces perhaps
the act of $A$ and the above splitting become of interest.

\hspace{0.25cm}
\section{Alternative (Oriented) singular homology}

It is not surprising that the existence of the alternative cohomology leads to an alternative homology. To see this
we must first construct alternative chains.  Usual singular chains form a free abelian group with generators all continuous
maps $\sigma:\Delta^{n}\xrightarrow{} X$ for all $n \geq 0$, namely for each $n \geq 0$, $C_n(X)=<\sigma:\Delta^{n}\xrightarrow{} X>$.
Of course for each $\sigma:\Delta^{n}\xrightarrow{} X$ according to Eq.~(\ref{perm2}), there are $(n+1)!$ maps  $\overline{s}^{\ast}\sigma$ for $s\in S_{n+1}$
and alternative chains are defined by inserting some relations on these maps as in the following definition.

\begin{definition} For each $n\geq 0$, alternative chains $C_{A  n}(X)$ are obtained from the group  $C_{n}(X)$
by considering the relations:

\begin{equation}\nonumber
\overline{s}^{\ast}\sigma=\epsilon(s)\sigma.
\end{equation}
for each $s\in S_{n+1}$, or briefly $C_{An}(X)=<\sigma:\Delta^{n}\xrightarrow{} X | \overline{s}^{\ast}\sigma=\epsilon(s)\sigma>$.

\end{definition}
According to notations of Ref. \cite{Barr}, $C_{A  n}(X)=C_n(X)/U_n(X)$ in which $U_n(X)$ is the subgroup of $C_n(X)$
generated by all expressions $(\overline{s}^{\ast}\sigma-\epsilon(s)\sigma)$ for all continuous maps $\sigma$ and all
permutations $s\in S_{n+1}$. Therefore each member of $C_{A  n}(X)$ is a coset with respect to $U_n(X)$ which we denote
it here by $[\sigma]\in C_{A  n}(X)=C_n(X)/U_n(X)$. Thus
\begin{equation}\label{relation}
[\overline{s}^{\ast}\sigma]=\epsilon(s)[\sigma].
\end{equation}

It was historically known that $C_{A  n}(X)$ is no longer a free group and contains torsion subgroups of order 2.
In fact if $\sigma:\Delta^{n}\xrightarrow{} X$ has some symmetry such that for an odd permutation $s$ ($\epsilon (s)=-1$)
we have $\overline{s}^{\ast}\sigma =\sigma$ then $2[\sigma ]=0$. Note that the situation does not change if these
conditions hold for one odd permutation or for more than one odd permutation and for both cases the class $[\sigma]=
[\overline{s}^{\ast}\sigma]$ for all $s\in S_{n+1}$ and has order 2. Of course if $\overline{s}^{\ast}\sigma =\sigma$
with $s$ an even permutation there is no torsion. To summarize, for each $n$, $C_{A  n}(X)$ splits into some free groups
and some torsion groups each isomorphic with $\mathbb{Z}_2$. Generators of the free part are all $[\sigma]$s such that
the symmetry $\overline{s}^{\ast}\sigma =\sigma$ either does not hold or $s$ is an even permutation and if $s$ is an odd
permutation, $[\sigma]$ is the generator of one $\mathbb{Z}_2$ summand.

Here in this section an outline is given for a proof of the isomorphism between the alternative and singular homologies.
Although alternative cochains are duals of alternative chains:
\begin{equation}\label{dual}
C^{n}_{A}(X;G)=(C_{An}(X))^{\ast}=\mathrm{Hom}(C_{A n}(X),G),
\end{equation}
unfortunately this isomorphism can not be extended to cohomology groups because alternative chains are not free abelian
groups. However, Barr demonstrated a homotopy equivalence between singular and alternative chain complexes which gives
the required isomorphisms \cite{Barr}. The next theorem gives a direct and simple demonstration
(different from those presented before) for the existence of the alternative (oriented) homology .

\begin{theorem} For a topological space $X$ the following sequence is a chain complex
\begin{equation}\nonumber
\dots\xrightarrow{\partial}C_{An}(X)\xrightarrow{\partial} C_{A(n-1)}(X)\xrightarrow{\partial} \dots\xrightarrow{\partial}  C_{A0}=C_{0}(X)\xrightarrow{\text{}}0.
\end{equation}
Therefore, the alternative homology $\mathrm{H}_{An}(X)=Ker\partial / Im \partial$ is defined.
\end{theorem}
{\it Proof}. Assume $\sigma:\Delta^{n}\xrightarrow{} X$ $(n\geq 1)$ is a representator out of $(n+1)!$ maps $\{\overline{s}^{\ast}\sigma \mid s\in S_{n+1}\}$ on which the relations $[\overline{s}^{\ast}\sigma ]=\epsilon(s)[\sigma ]$ are imposed. We must show that
for each $s' \in S_{n}$, $[\overline{s'}^{\ast}(\partial \sigma)]=\epsilon(s')[\partial\sigma]$ which indicates that the boundary operator $\partial$ is also well defined on $C_{An}(X)$ with the image in $C_{A(n-1)}(X)$. For $s'\in S_{n}$ ,

\begin{equation}\label{12}
\overline{s'}^{\ast}(\partial \sigma)
=\sum_{i=0}^{n}(-1)^{i}(\overline{s'}^{\ast}\sigma|_{[v_{0},\ldots,\hat{v}_{i},\ldots,v_{n}]}).
\end{equation}

In the proof of Theorem 4, permutations $s_{}^{(i)}(j,s')\in S_{n+1}$ were defined such that $s_{i}^{(i)}(j,s')=s'$ and $(s^{(i)}(j,s'))(i)=j$.
Now especially put $j=i$ and consider $s^{(i)}(i,s')$ which for an easier notation we denote it by $s^{(i)}$. Then $s_{i}^{(i)}=s'$
($s_{i}$ is defined by Definition 2), $s^{(i)}_{}(i)=i$ and Eq.~(\ref{12}) becomes

\begin{equation}\nonumber
\begin{split}
\overline{s'}^{\ast}(\partial\sigma)=\sum_{i=0}^{n}(-1)^{i}\overline{s_{i}^{(i)}}^{\ast}(\sigma|_{[v_{0},\ldots,\hat{v}_{s^{(i)}(i)},\ldots,v_{n}]})=\sum_{i=0}^{n}(-1)^{i}(\overline{s^{(i)}}^{\ast}\sigma)|_{[v_{0},\ldots,\hat{v}_{i},\ldots,v_{n}]}.
\end{split}
\end{equation}
Since  $[\overline{s^{(i)}}^{\ast}\sigma]=\epsilon(s^{(i)})[\sigma]$  and $\epsilon(s^{(i)})=(-1)^{i-i}\epsilon(s')=\epsilon(s')$  (Eq.~(\ref{sigma})), we obtain

\begin{equation}\nonumber
\begin{split}
[(\overline{s'}^{\ast})(\partial\sigma)]=\epsilon(s')[\sum_{i=0}^{n}(-1)^{i}\sigma|_{[v_{0},\ldots,\hat{v}_{i},\ldots,v_{n}]}]=
\epsilon(s')[\partial\sigma].
\end{split}
\end{equation} $\blacksquare$\\
It is also seen that
\begin{equation}\label{hom0}
\mathrm{H}_{A0}(X)=\mathrm{H}_{0}(X).
\end{equation}

Although the above theorem is complete but it is useful to see directly that the boundary operator $\partial$ as a homomorphism
takes each torsion subgroup to a torsion subgroup and not to a free subgroup of the successive alternative chain group. Suppose
$[\sigma]\in C_{A n}(X)$ is the generator of a $\mathbb{Z}_2$ summand and thus there is an odd permutation $s\in S_{n+1}$ such
that $\overline{s}^{\ast}\sigma =\sigma$, or equivalently $\sigma =-\epsilon(s)\overline{s}^{\ast}\sigma$. Taking the boundary
of this equation and using Eq.~(\ref{5}) we find
\begin{equation}\nonumber
\partial \sigma =-\sum_{i=0}^{n}(-1)^{i}\epsilon(s)\overline{s}_{i}^{\ast}(\sigma|_{[v_{0},\ldots,\hat{v}_{s(i)},\ldots,v_{n}]}).
\end{equation}
We then substitute for $\epsilon(s)$ from Eq.~(\ref{sigma}) to obtain
\begin{equation}\nonumber
\partial \sigma =-\sum_{i=0}^{n}(-1)^{s(i)}\epsilon(s_i)\overline{s}_{i}^{\ast}(\sigma|_{[v_{0},\ldots,\hat{v}_{s(i)},\ldots,v_{n}]}).
\end{equation}
We take the coset from the both sides of the above equation regarding the fact that due to Relation (\ref{relation}) we have
$\epsilon(s_i)[\overline{s}_{i}^{\ast}(\sigma|_{[v_{0},\ldots,\hat{v}_{s(i)},\ldots,v_{n}]})]=[\sigma|_{[v_{0},\ldots,\hat{v}_{s(i)},\ldots,v_{n}]}]$
and finally obtain
\begin{equation}\label{bound-torsion1}
[\partial \sigma ]=-\sum_{i=0}^{n}(-1)^{s(i)}[\sigma|_{[v_{0},\ldots,\hat{v}_{s(i)},\ldots,v_{n}]}]=-[\partial \sigma],
\end{equation}
and thus $[\partial \sigma]$ is a torsion member of $C_{A (n-1)}(X)$. Even more precisely we can compare the coset of each individual term
in $\partial \sigma$: Since $\sigma =\overline{s}^{\ast}\sigma$ so
\begin{equation}\nonumber
\sigma|_{[v_{0},\ldots,\hat{v}_i,\ldots,v_{n}]}=(\overline{s}^{\ast}\sigma )|_{[v_{0},\ldots,\hat{v}_i,\ldots,v_{n}]}=
\overline{s}_{i}^{\ast}(\sigma|_{[v_{0},\ldots,\hat{v}_{s(i)},\ldots,v_{n}]}),
\end{equation}
which taking their cosets with the use of (\ref{relation}) and (\ref{sigma}) (where $\epsilon(s)=-1$) yields
\begin{equation}\label{bound-torsion2}
(-1)^i[\sigma|_{[v_{0},\ldots,\hat{v}_i,\ldots,v_{n}]}]=(-1)^i\epsilon(s_i)[\sigma|_{[v_{0},\ldots,\hat{v}_{s(i)},\ldots,v_{n}]}]
=-(-1)^{s(i)}[\sigma|_{[v_{0},\ldots,\hat{v}_{s(i)},\ldots,v_{n}]}].
\end{equation}
This means that the $i$-th term in the expression for $[\partial \sigma]$ is canceled with the $s(i)$-th term in that expression. Thus all
pairs $i$ and $s(i)$ with $s(i)\neq i$ cancel each other and only remain terms with $s(i)=i$ which by the above equation are of order 2. Hence
taking the boundary of a generator of order 2 gives a summation of generators each with order 2.

For any continuous map $f:X\xrightarrow{}Y$, $f_{\ast}:\mathrm{H}_{An}(X)\xrightarrow{}\mathrm{H}_{An}(Y)$ is also well defined. Also all exact sequences such as for a pair $(X,B)$ $(B\subset X)$ or
for a triple $(C\subset B\subset X)$ are valid for alternative homologies too. It is also seen that the  excision isomorphism holds.
Indeed this isomorphism is proved based on a barycentric subdivision (see for example Ref. \cite{Hatcher}) which is derived inductively from faces of $\Delta^{n}$.
In the proof of Theorem 5 we obtained that
\begin{equation}\label{14}
[\overline{s'}^{\ast}(\sigma|_{F_i})]=[(\overline{s^{(i)}}^{\ast}\sigma)|_{F_i}]=\epsilon(s^{(i)})[\sigma |_{F_i}]=\epsilon(s')[\sigma |_{F_i}],
\end{equation}
which shows that the restriction of $\sigma$ to any face $F_{i}$ also satisfies the relations necessary to construct $C_{A(n-1)}(X)$.
A detailed and precise review of the excision theorem admits to hold it also for  alternative chains due to Eq.~(\ref{14}). Excision theorem
leads to the Mayer-Vietoris  exact sequence for alternative homologies.

Now if $F:X \times I \xrightarrow{}Y$ be a homotopy between $f , g: X\xrightarrow{}Y$, there is the famous chain map
$P:C_{n}(X)\xrightarrow{}C_{n+1}(Y)$ defined as
\begin{equation}\label{prism}
P(\sigma)=\sum_{i=0}^{n}(-1)^i(F\circ(\sigma\times id))|_{[v_0,\ldots,v_i,w_i,\ldots,w_n]},
\end{equation}
where $v_0,\ldots,v_n$ are the vertices of $\Delta^n\times\{0\}$ and $w_0,\ldots,w_n$ are the vertices of $\Delta^n\times\{1\}$ as
subsets of $\Delta^n\times I$. It then can be shown that $\partial P+ P \partial=g_{\#}-f_{\#}$ \cite{Hatcher}. In order to modify
this operator to be well defined between alternative chains we extend the function $F$ to $\tilde{F}$ defined on chains such that
for each $n$ and $s\in S_{n+1}$  we have
\begin{equation}\label{Ftilde}
\begin{split}
(\tilde{F}\circ(\overline{s}^{\ast}\sigma\times id))|_{[v_0,\ldots,v_i,w_i,\ldots,w_n]}=(\tilde{F}\circ(\sigma\times id))|_{[v_{s(0)},\ldots,v_{s(i)},w_{s(i)},\ldots,w_{s(n)}]}\\
=\epsilon (s)(\tilde{F}\circ(\sigma\times id))|_{[v_0,\ldots,v_i,w_i,\ldots,w_n]}.
\end{split}
\end{equation}
Now we take the cosets and finally define $\tilde{P}:C_{An}(X) \xrightarrow{}C_{A(n+1)}(Y)$ as
\begin{equation}\label{prismtilde}
\tilde{P}([\sigma ])=\sum_{i=0}^{n}(-1)^i[(\tilde{F}\circ(\sigma\times id))|_{[v_0,\ldots,v_i,w_i,\ldots,w_n]}].
\end{equation}
It is not difficult to observe that taking the cosets in the above equation is consistent with definition (\ref{Ftilde}). Indeed, 
assume the permutation imposed on $n+2$ vertices $[v_0,\ldots,v_i,w_i,\ldots,w_n]$ is such that $v_i$ and $w_i$ remain unchanged and 
there are only permutations between $v_0,\ldots,v_{i-1}$ from one hand and between $w_{i+1},\ldots,w_n$ form the other hand. Then this
permutation is practically in $S_{n+1}$ defined in (\ref{Ftilde}) and so $\epsilon (s)$ is in agreement with this definition 
and the definition of cosets (\ref{relation}). If the permutation is not like the above then it is different from those considered
in (\ref{Ftilde}) and so there is no conflict between this definition and (\ref{relation}). Hence, the chain map (\ref{prismtilde})
is well defined and of course satisfy the relation $\partial\tilde{P}+ \tilde{P}\partial=g_{\#}-f_{\#}$ on alternative chain complexes.
Hence, $f_{\ast}=g_{\ast}$ on $\mathrm{H}_{An}(X)$ and since alternative cochains are dual to alternative chains, $f^{\ast}=g^{\ast}$ on $\mathrm{H}_{A}^{n}(X)$ as well. This duality also establishes the Mayer-Vietoris long exact sequence for alternative cohomologies 
provided that this sequence holds for alternative homologies. We recall how to obtain
homology groups of an $n$-sphere $S^n:\mathrm{H}_{0}(S^n)=\mathbb{Z}$ and for a good pair $(D^n,\partial D^n=S^{{n-1}})$ we use the
related long exact sequence and the fact that $\mathrm{H}_{k}(D^n,\partial D^n)=\mathrm{H}_{k}(D^n/\partial D^n)=\mathrm{H}_{k}(S^n)$.
Now since $\mathrm{H}_{A0}(X)=\mathrm{H}_{0}(X)$ and due to $f_\ast=g_\ast$ (for $f\thicksim g$) we find that

\begin{equation}\label{15}
\begin{split}
\mathrm{H}_{Ak}(S^n)=\mathrm{H}_{k}(S^n),
\end{split}
\end{equation}
 for all $k$ and $n$.

 A $\Delta$-complex structure on $X$ is a collection of characteristic maps $\{\sigma_{\alpha}:\Delta^{n_{\alpha}}\xrightarrow{}X\}$
 such that for each $\alpha$ and each $i$, $\sigma_{\alpha} |_{F_i}=\sigma_{\beta}$ for some $\beta$ ($n_{\beta}=n_{\alpha}-1$)
 and $\sigma_{\alpha}$ is a homeomorphism when restricted to the interior of $\Delta^{n_{\alpha}}$. If in addition for each $\alpha$
 the images of the vertices of $\Delta^{n_{\alpha}}$ under $\sigma_{\alpha}$ are unique for $\sigma_{\alpha}$, i.e. no other map $\sigma_{\beta}$ exists such that $\{\sigma_{\alpha}(v_i)|i=0,\dots, n_{\alpha}\}=\{\sigma_{\beta}(v_i)|i=0,\dots, n_{\beta}=n_{\alpha}\}$, then this $\Delta$-complex
 is called a simplicial complex. A simplicial homology is defined for a $\Delta$-complex or more specially for a simplicial complex.
 In any text book such as Ref. \cite{Hatcher} the isomorphism between the singular and simplicial homologies is proved.
 This proof uses the exact sequence for pair $(\Delta^{n},\partial\Delta^{n})$ and the homology groups of spheres and the five-lemma
 which all of these tools are valid for alternative homologies too. Thus we see that for any simplicial complex structure on $X$, the
 alternative homology is isomorphic to the usual simplicial homology which is itself isomorphic to the singular homology. Here we obtain:
\begin{theorem}
For a simplicial complex structure on $X$ the alternative and singular homologies are isomorphic.
\end{theorem}
On the other hand any $CW$ structure is homotopy equivalent to a simplicial complex \cite{Hatcher}. Thus
\begin{corollary} For any $CW$ complex, singular and alternative homologies are isomorphic.
\end{corollary}

Finally any topological space $X$ is weakly homotopy equivalent to a $CW$ space which induces isomorphisms
on homology groups \cite{Hatcher}. These isomorphisms also hold for alternative homology groups and so

\begin{corollary} For any topological space $X$,

\begin{equation}\nonumber
\mathrm{H}_{An}(X)=\mathrm{H}_{n}(X).
\end{equation}
\end{corollary}

As mentioned before, the isomorphisms of cohomology groups can not follow here since alternative chain groups are not free
and Ref. \cite{Barr} does this job. Again we emphasize that a detailed proof of the above results requires to add a long
unnecessary well known methods in text books while regarding the comments here for alternative chains is enough for a reader
familiar to those standard text book methods.

 \hspace{0.25cm}
\section{Summary and conclusive remarks}

The main motivative of this paper was to get close to the answer to this question that to what extent a nonlinear differential equation
does have solutions independent from the given structure to a manifold? Then it was tried to create a duality between a nonlinear
differential equation on a manifold and a nonlinear cochain equation based only on topological properties. Any correspondence between
the solutions of these two equations may be the subject of future investigations such as equation $\alpha \wedge d\alpha=0$ investigated
by the authors and will appear soon.

Similar to definition of exterior (alternative) forms in linear algebra, alternative cochains as a subgroup of singular cochains were 
considered which historically known as oriented singular cochains. If the coefficients of these cochains are in rational or real numbers, 
then an alternative maker operator (homomorphism) naturally appeared by which a splitting of singular cochains was induced whose alternative 
cochains was one of its summands. This alternative maker operator was also used to define a modified cup product only for alternative 
cochains very similar to the wedge product in differential forms. This alternative cup product was seen to have all three algebraic properties 
holding for the wedge product.

Although it was known that these alternative cochains of all orders construct a chain complex by the usual coboundary operator, here a new
and more direct proof for this fact were presented. So alternative cohomologies were defined. Again if the coefficients are rational or 
real numbers, the alternative maker operator (prove here to commute with the coboundary operator) holds also for cohomologies resulting 
a splitting of singular cohomology groups to alternative cohomologies and another summand. Since it was already proved that alternative 
and singular cohomologies are isomorphic, the second summand in the above splitting is zero if the topological space is compact. If the 
considered space is not compact, maybe the second summand is not zero and can be a subject for further investigations. 

The existence of alternative chain complex (known historically as oriented singular chain complex) were explicitly shown in a new and direct way. 
It was shown directly that all long exact sequences and homotopy properties or briefly, all homology and cohomology properties hold for alternative 
homology and cohomology as well. Finally the isomorphism between alternative and singular homologies was justified without going to long details. 
However this approach is not able to give the isomorphism between the corresponding cohomologies and this task was already done by Barr by proving that the singular and oriented (alternative) chain complexes are chain homotopy equivalent \cite{Barr}.

Therefore alternative cochains possessing algebraic properties similar to differential forms give alternative cohomology which is the same
as singular cohomology. The isomorphism between the de Rham and singular cohomologies suggests to investigate future dualities for other
differential equations probably nonlinear.

The last subject of interest is to seek for a generalized characteristic class of a vector bundle: Any vector bundle selects special solutions
of $\delta\alpha=0$, out of the members of the cohomology groups. Does this bundle select some solutions of other nonlinear cochain
equations such as $\alpha\overset{A}\smile\delta\alpha=0$?

 \hspace{0.25cm}
\section{Acknowledgments}

The second author would like to appreciate the partial financial support and kind hospitality of School of Mathematics, Institute for Research in Fundamental Sciences (IPM). The authors wishes to acknowledge Eaman Eftekhary, Ali Kamalinejad and Iman Setayesh for their remarks.

\hspace{0.25cm}

\end{document}